\documentclass{amsart}[10pt]
\usepackage{a4}
\usepackage[T1]{fontenc}
\usepackage[latin1]{inputenc}
\usepackage{hyperref}
\usepackage[dvips]{graphics}
\usepackage[english]{babel}
\usepackage{amsmath}
\usepackage{amssymb}
\usepackage{amsopn}
\usepackage{amsbsy}
\usepackage{amstext}
\usepackage{latexsym}
\usepackage{amsfonts}
\usepackage{array}
\usepackage{epsfig}
\usepackage{mathrsfs}
\title[A curious result related to Kempner's series]{A curious result related to Kempner's series}
\author[BAKIR FARHI]{Bakir FARHI}

\newtheorem{thm}{Theorem}
\newtheorem{prop}[thm]{Proposition}

\let\epsilon=\varepsilon

\def\EMdash{\leavevmode\hbox to 7.5mm{\vrule height .63ex depth -.59ex
    width 5.4mm\hfill}}

\begin{document}
\maketitle \vspace{-4cm}
\begin{flushleft}
To appear in {\it Amer. Math. Monthly}.
\end{flushleft}~\vspace{2cm}

\begin{center}
Département de Mathématiques, Université du Maine, \\
Avenue Olivier Messiaen, 72085 Le Mans Cedex 9, France. \\
bakir.farhi@gmail.com
\end{center}
\vskip 1cm \noindent {\small {\bf Abstract.---} It is well known
since A. J. Kempner's work that the series of the reciprocals of
the positive integers whose the decimal representation does not
contain any digit $9$, is convergent. This result was extended by
F. Irwin and others to deal with the series of the reciprocals of
the positive integers whose the decimal representation contains
only a limited quantity of each digit of a given nonempty set of
digits. Actually, such series are known to be all convergent.
Here, letting $S^{(r)}$ $(r \in \mathbb{N})$ denote the series of
the reciprocal of the positive integers whose the decimal
representation contains the digit $9$ exactly $r$ times, the
impressive obtained result is that $S^{(r)}$ tends to $10
\log{10}$ as $r$ tends to infinity!} \vskip 2mm \noindent
\hfill{--------} 

\section{Introduction}

Throughout this article, we let $\mathbb{N}^*$ denote the set
$\mathbb{N} \setminus \{0\}$ of positive integers. We let
$E^{(r)}$ $(r \in \mathbb{N})$ denote the set of all positive
integers whose decimal representation contains the digit $9$
exactly $r$ times. We also let $E_{n}^{(r)}$ $(r , n \in
\mathbb{N})$ denote the set $E^{(r)} \cap [10^n , 10^{n + 1}[$. We
clearly have for any $r \in \mathbb{N}$:
$$E^{(r)} = \bigsqcup_{n \in \mathbb{N}} E_{n}^{(r)} .$$
For all $r , n \in \mathbb{N}$, write:
$$S_{n}^{(r)} := \sum_{k \in E_{n}^{(r)}} \frac{1}{k}$$
and for all $r \in \mathbb{N}$, write:
$$S^{(r)} := \sum_{k \in E^{(r)}} \frac{1}{k} = \sum_{n \in \mathbb{N}} S_{n}^{(r)} .$$

In 1914, A. J. Kempner \cite{kem} showed that the series $S^{(0)}$
converges. After him, several generalizations were obtained by
several authors. Among others, F. Irwin \cite{irw} showed that the
series derived from the harmonic series $1 + \frac{1}{2} +
\frac{1}{3} + \cdots$ by including only those terms whose
denominators contain a limited quantity of each digit of a given
nonempty set of digits is convergent. It follows in particular
that the series $S^{(r)}$ converges for all natural numbers $r$.
In this paper, the main result obtained is that the sequence
${(S^{(r)})}_{r \geq 1}$ decreases and converges to $10 \log{10}$
(see Theorem \ref{t2}). As a consequence, we deduce that we have
$S^{(r)} > 10 \log{10} \simeq 23.025$ $(\forall r \geq 1)$. So,
according to the calculations of R. Baillie \cite{bai}, we have
the unexpected inequality $S^{(0)} \simeq 22.920 < S^{(1)}$. We
must notice that the approximate numerical values of the
$S^{(r)}$'s $(r \geq 1)$ are very difficult to calculate.

In the last section of the paper, we state a generalization of our
main result by taking instead of the digit $9$ any other digit $d
\in \{0 , 1 , \dots , 9\}$.

\section{The Results}

Suppose $r , n \in \mathbb{N}^*$. If a positive integer $k$
belongs to $E_{n}^{(r)}$, then writing
$$k = 10 t + \ell$$
with $t \in \mathbb{N}$ and $\ell \in \{0 , 1 , \dots , 9\}$, we
clearly have either

$\bullet$ $t \in E_{n - 1}^{(r)}$ and $\ell \in \{0 , 1 , \dots ,
8\}$, or

$\bullet$ $t \in E_{n - 1}^{(r - 1)}$ and $\ell = 9$.

\noindent It follows that
\begin{equation}\label{eq1}
S_{n}^{(r)} = \sum_{\ell = 0}^{8} \sum_{t \in E_{n - 1}^{(r)}}
\frac{1}{10 t + \ell} + \sum_{t \in E_{n - 1}^{(r - 1)}}
\frac{1}{10 t + 9} .
 \end{equation}
 We will find it useful to approximate $S_{n}^{(r)}$ with the
 simpler formula
\begin{equation}\label{eq2}
T_{n}^{(r)} := \sum_{\ell = 0}^{8} \sum_{t \in E_{n - 1}^{(r)}}
\frac{1}{10 t} + \sum_{t \in E_{n - 1}^{(r - 1)}} \frac{1}{10 t} =
\frac{9}{10} S_{n - 1}^{(r)} + \frac{1}{10} S_{n - 1}^{(r - 1)} .
\end{equation}
The error in this approximation is given by
\begin{equation}\label{eq5}
C_{n , r} := T_n^{(r)} - S_{n}^{(r)} .
\end{equation}
So, we have
\begin{equation}\label{eq12}
S_{n}^{(r)} = \frac{9}{10} S_{n - 1}^{(r)} + \frac{1}{10} S_{n -
1}^{(r - 1)} - C_{n , r} .
\end{equation}
This identity will play an important role in what follows.

Our first proposition shows that the errors $C_{n , r}$ are not
very large.
\begin{prop}\label{t1}
 The real numbers $C_{n , r}$ $(r , n \in \mathbb{N}^*)$  are all
 nonnegative and we have
$$\sum_{r = 1}^{\infty} \sum_{n = 1}^{\infty} C_{n , r} < \infty .$$
\end{prop}\noindent
{\bf Proof.} From (\ref{eq1}) and the first equality of
(\ref{eq2}), we have for all $r , n \in \mathbb{N}^*$:
\begin{equation}\label{eq3}
C_{n , r} = \sum_{\ell = 1}^{8} \sum_{t \in E_{n - 1}^{(r)}}
\frac{\ell}{10 t (10 t + \ell)} + \sum_{t \in E_{n - 1}^{(r - 1)}}
\frac{9}{10 t (10 t + 9)} .
\end{equation}
This last identity shows that we have $C_{n , r} \geq 0$ $(\forall
r , n \in \mathbb{N}^*)$. Moreover, using (\ref{eq3}), we have for
all $r , n \in \mathbb{N}^*$:
$$
C_{n , r} \leq \frac{9}{25} \sum_{t \in E_{n - 1}^{(r)}}
\frac{1}{t^2} + \frac{9}{100} \sum_{t \in E_{n - 1}^{(r - 1)}}
\frac{1}{t^2} .
$$
Since the sets $E_{n}^{(r)}$ $(r , n \in \mathbb{N})$ form a
partition of $\mathbb{N}^*$, it follows that:
$$\sum_{r = 1}^{\infty} \sum_{n = 1}^{\infty} C_{n , r} \leq \left(\frac{9}{25} + \frac{9}{100}\right)
\sum_{t = 1}^{\infty} \frac{1}{t^2} = \frac{9}{20} \frac{\pi^2}{6}
< \infty .$$ The proof is
complete.\penalty-20\null\hfill$\blacksquare$\par\medbreak

For the following, put for all $r \in \mathbb{N}^*$:
$$C_r := \sum_{n = 1}^{\infty} C_{n , r} .$$
According to (\ref{eq3}), we clearly have $C_{n , r} > 0$ for all
$n , r \in \mathbb{N}^*$ such that $n \geq r$. Consequently, we
have $C_r > 0$ for all $r \in \mathbb{N}^*$.
\begin{prop}\label{coll1}
For all $r \in \mathbb{N}$, the series $S^{(r)}$ converges. In
addition, we have:
\begin{equation}\label{eq6}
S^{(1)} = S^{(0)} - 10 C_1 + \frac{10}{9}
\end{equation}
and
\begin{equation}\label{eq7}
S^{(r)} = S^{(r - 1)} - 10 C_r ~~~~~~ (\forall r \geq 2) .
\end{equation}
\end{prop}\noindent
{\bf Proof.} The fact that the series $S^{(r)}$ $(r \in
\mathbb{N})$ are all convergent is already known (see, e.g.,
\cite{irw}). Let us prove the relations (\ref{eq6}) and
(\ref{eq7}) of the proposition. Using the relation (\ref{eq12}),
we have for all $r \in \mathbb{N}^*$:
\begin{eqnarray*}
S^{(r)} & = & \sum_{n = 1}^{\infty} S_n^{(r)} + S_0^{(r)} \\
& = & \sum_{n = 1}^{\infty}\left(\frac{9}{10} S_{n - 1}^{(r)} +
\frac{1}{10} S_{n - 1}^{(r - 1)} - C_{n , r}\right) + S_0^{(r)} \\
& = & \frac{9}{10} \sum_{n = 1}^{\infty} S_{n - 1}^{(r)} +
\frac{1}{10} \sum_{n = 1}^{\infty} S_{n - 1}^{(r - 1)} - \sum_{n =
1}^{\infty} C_{n , r} + S_0^{(r)} \\
& = & \frac{9}{10} \sum_{n = 0}^{\infty} S_n^{(r)} + \frac{1}{10}
\sum_{n = 0}^{\infty} S_n^{(r - 1)} - \sum_{n =
1}^{\infty} C_{n , r} + S_0^{(r)} \\
& = & \frac{9}{10} S^{(r)} + \frac{1}{10} S^{(r - 1)} - C_r +
S_0^{(r)} .
\end{eqnarray*}
Hence:
\begin{eqnarray*}
S^{(r)} & = & S^{(r - 1)} - 10 C_r + 10 S_0^{(r)} \\
& = & \begin{cases} S^{(r - 1)} - 10 C_r + \frac{10}{9} & \text{if
$r = 1$}, \\ S^{(r - 1)} - 10 C_r & \text{if $r \geq
2$}.\end{cases}
\end{eqnarray*}
This confirms the required formulas (\ref{eq6}) and (\ref{eq7}) of
the proposition and finishes this
proof.\penalty-20\null\hfill$\blacksquare$\par\medbreak

We now arrive at the most important and completely new result of
this paper:
\begin{thm}\label{t2}
The sequence ${(S^{(r)})}_{r \geq 1}$ decreases and converges to
$10 \log{10}$. In particular, we have:
$$S^{(r)} > 10 \log{10} ~~~~~~ (\forall r \geq 1) .$$
\end{thm}\noindent
{\bf Proof.} Since $C_r > 0$ $(\forall r \geq 1)$, the formula
(\ref{eq7}) of Proposition \ref{coll1} shows that the sequence
${(S^{(r)})}_{r \geq 1}$ decreases. Since this sequence is
positive (so bounded from below by $0$), it is necessarily
convergent. It remains to calculate its limit as $r$ tends to
infinity. We have for all integer $R \geq 2$:
\begin{eqnarray*}
S^{(R)} & = & \sum_{r = 2}^{R} \left(S^{(r)} - S^{(r - 1)}\right)
+ S^{(1)} \\
& = & \sum_{r = 2}^{R} (- 10 C_r) + S^{(0)} - 10 C_1 +
\frac{10}{9} ~~~~~~ \text{(according to (\ref{eq6}) and
(\ref{eq7}))} \\
& = & S^{(0)} - 10 \sum_{r = 1}^{R} C_r + \frac{10}{9} .
\end{eqnarray*}
Hence:
\begin{equation}\label{eq8}
\lim_{R \rightarrow \infty} S^{(R)} = S^{(0)} - 10 \sum_{r = 1}^{
\infty} C_r + \frac{10}{9} .
\end{equation}

Now, let us calculate the sum $\sum_{r = 1}^{\infty} C_r = \sum_{r
= 1}^{\infty} \sum_{n = 1}^{\infty} C_{n , r}$. From (\ref{eq1}),
(\ref{eq2}), and (\ref{eq5}), we have for all $r , n \in
\mathbb{N}^*$:
\begin{eqnarray*}
C_{n , r} & := & T_{n}^{(r)} - S_{n}^{(r)} \\
& = & \sum_{\ell = 0}^{8} \sum_{t \in E_{n - 1}^{(r)}}
\left(\frac{1}{10 t} - \frac{1}{10 t + \ell}\right) + \sum_{t \in
E_{n - 1}^{(r - 1)}} \left(\frac{1}{10 t} - \frac{1}{10 t +
9}\right) .
\end{eqnarray*}
By remarking that the sets $E_{n - 1}^{(r)}$ $(r , n \in
\mathbb{N}^*)$ form a partition of $\mathbb{N}^* \setminus
E^{(0)}$ and that the sets $E_{n - 1}^{(r - 1)}$ $(r , n \in
\mathbb{N}^*)$ form a partition of $\mathbb{N}^*$, we deduce that:
$$\sum_{r = 1}^{\infty} \sum_{n = 1}^{\infty} C_{n , r} ~~=~~ \sum_{\ell = 0}^{8} \sum_{t \in \mathbb{N}^*
\setminus E^{(0)}} \left(\frac{1}{10 t} - \frac{1}{10 t +
\ell}\right) + \sum_{t = 1}^{\infty} \left(\frac{1}{10 t} -
\frac{1}{10 t + 9}\right)$$
\begin{equation*}
\begin{split}
&= \sum_{\ell = 0}^{8}\left\{\sum_{t = 1}^{
\infty}\left(\frac{1}{10 t} - \frac{1}{10 t + \ell}\right) -
\sum_{t \in E^{(0)}}\left(\frac{1}{10 t} - \frac{1}{10 t +
\ell}\right)\right\} + \sum_{t = 1}^{\infty}\left(\frac{1}{10 t}
- \frac{1}{10 t + 9}\right) \\
&= \sum_{\ell = 0}^{9}\sum_{t = 1}^{\infty}\left(\frac{1}{10 t} -
\frac{1}{10 t + \ell}\right) - \sum_{\ell = 0}^{8}\sum_{t \in
E^{(0)}}\left(\frac{1}{10 t} - \frac{1}{10 t + \ell}\right) \\
&= \sum_{\ell = 0}^{9}\sum_{t = 1}^{\infty}\left(\frac{1}{10 t} -
\frac{1}{10 t + \ell}\right) - \frac{9}{10}\sum_{t \in E^{(0)}}
\frac{1}{t} + \sum_{\ell = 0}^{8}\sum_{t \in E^{(0)}} \frac{1}{10
t + \ell} .
\end{split}
\end{equation*}
Since $\sum_{t \in E^{(0)}} 1/t = S^{(0)}$ and
$$\sum_{\ell = 0}^{8}\sum_{t \in E^{(0)}} \frac{1}{10 t + \ell} =
\sum_{t \in E^{(0)}} \frac{1}{t} - \sum_{t = 1}^{8} \frac{1}{t} =
S^{(0)} - \sum_{t = 1}^{8} \frac{1}{t}$$ (because $\{10 t + \ell
~|~ t \in E^{(0)} , \ell \in \{0 , 1 , \dots , 8\}\} = E^{(0)}
\setminus \{1 , 2 , \dots , 8\}$), it follows that:
\begin{equation}\label{eq9}
\sum_{r = 1}^{\infty}\sum_{n = 1}^{\infty} C_{n , r} = \Delta +
\frac{1}{10} S^{(0)} - \sum_{t = 1}^{8} \frac{1}{t} ,
\end{equation}
where
$$
\Delta := \sum_{\ell = 0}^{9}\sum_{t = 1}^{
\infty}\left(\frac{1}{10 t} - \frac{1}{10 t + \ell}\right) =
\sum_{t = 1}^{\infty}\left\{\frac{1}{t} - \sum_{10 t \leq m < 10
(t + 1)} \frac{1}{m}\right\} .$$

Let us calculate $\Delta$. For all sufficiently large positive
integers $N$, we have:
$$
\sum_{t = 1}^{N}\left\{\frac{1}{t} - \sum_{10 t \leq m < 10 (t +
1)} \frac{1}{m}\right\} ~=~ \sum_{t = 1}^{N} \frac{1}{t} - \sum_{m
= 10}^{10 N + 9} \frac{1}{m} \phantom{aaaaaaaaaaaaaaaaaaaaaaaaaaa}
$$
\begin{equation*}
\begin{split}
&= \left(\log{N} + \gamma\right) - \left\{\log(10 N + 9) + \gamma
- \sum_{m = 1}^{9} \frac{1}{m}\right\} + o_N(1) \\
&= \sum_{m = 1}^{9} \frac{1}{m} - \log{10} + o_N(1) ,
\end{split}
\end{equation*}
where $\gamma$ denotes the Euler's constant. By taking the limits
as $N$ tends to infinity, we obtain:
$$\Delta = \sum_{m = 1}^{9} \frac{1}{m} - \log{10} .$$

Now, by substituting this value of $\Delta$ into (\ref{eq9}), we
obtain:
$$
\sum_{r = 1}^{\infty} C_r = \sum_{r = 1}^{\infty}\sum_{n = 1}^{
\infty} C_{n , r} = \frac{1}{10} S^{(0)} - \log{10} + \frac{1}{9}
.
$$
Finally, by substituting this value of $\sum_{r = 1}^{\infty} C_r$
into (\ref{eq8}), we conclude that:
$$\lim_{R \rightarrow \infty} S^{(R)} = 10 \log{10}$$
as required. The proof is
complete.\penalty-20\null\hfill$\blacksquare$\par\medbreak\noindent
\section{Generalization to other digits}
The method presented above can be applied to obtain the general
result related to any digit $d \in \{0 , 1 , \dots , 9\}$. For $d
\in \{0 , 1 , \dots , 9\}$ and for $r \in \mathbb{N}$, let
$\sigma_d(r)$ denote the sum of the reciprocals of all positive
integers whose decimal representation contains the digit $d$
exactly $r$ times. Then, we have the following:
\begin{thm}\label{t3}
For all $d \in \{0 , 1 , \dots , 9\}$, the sequence
$(\sigma_d(r))_{r \geq 1}$ decreases and converges to $10
\log{10}$.
\end{thm}\noindent
{\bf Remark.} According to Theorem \ref{t3} and to the approximate
values of the numbers $\sigma_d(0)$ $(0 \leq d \leq 9)$ given by
R. Baillie \cite{bai}, we see that only in the particular case $d
= 0$ (already studied by A. D. Wadhwa \cite{wad}) does the
sequence $(\sigma_d(r))_r$ start its decrease at $r = 0$.

\end{document}